\documentclass[12pt]{article}

\usepackage{amsfonts}
\usepackage{amsxtra}
\usepackage{amsthm}
\usepackage{amssymb}
\usepackage{amsmath,amscd}
\usepackage{supertabular}
\usepackage{pstricks}
\usepackage{pst-plot}
\usepackage{pstricks-add}
\usepackage{pst-eps}
\usepackage{makeidx}
\usepackage{bm}

\newcommand{\vect}[1]{\mathbf{#1}}
\newcommand{\abs}[1]{{\mathopen\mid}{#1}{\mathclose\mid}}

\newtheorem{theorem}{Theorem}[section]
\newtheorem{lemma}[theorem]{Lemma}
\newtheorem{definition}[theorem]{Definition}

\newtheorem{corollary}[theorem]{Corollary}
\newtheorem{question}[theorem]{Question}

\newtheorem{example}[theorem]{Example}

\DeclareMathOperator\GL{GL}
\DeclareMathOperator\PGL{PGL}

\DeclareMathOperator\PU{PU}

\DeclareMathOperator\Aut{Aut}
\DeclareMathOperator\Art{Art}
\DeclareMathOperator\Hom{Hom}
\DeclareMathOperator\Ker{Ker}

\newcommand{\R}{{\mathbb R}}
\newcommand{\N}{{\mathbb N}}
\newcommand{\C}{{\mathbb C}}
\newcommand{\Q}{{\mathbb Q}}
\newcommand{\Pro}{{\mathbb P}}

\newcommand{\Z}{{\mathbb Z}}
\newcommand{\B}{{\mathbb B}}

\newcommand{\UH}{{\mathbb H}}

\begin{document}

\title{Hyperbolic Structures and Root Systems}

\author
{Gert Heckman and Eduard Looijenga}

\maketitle

\begin{abstract}
We discuss the construction of a one parameter family of complex hyperbolic structures 
on the complement of a toric mirror arrangement associated with a simply laced root system.
Subsequently we find conditions for which parameter values this leads to ball quotients.
\end{abstract}

\section{Geometric structures}

Suppose $M$ is a connected smooth complex manifold of dimension $n$.
In fact we shall assume that $M$ is an open subset of $\C^n$ with coordinates $z=(z_1,\cdots,z_n)$ on $M$.
Write $\partial_i$ for $\partial/\partial z_i$.
A connection on the one jet bundle of $M$ is given by a system of $n(n+1)/2$ linear differential equations
\[ [\partial_i\partial_j+\sum_k\;a_{ij}^k(z)\partial_k+b_{ij}(z)]f(z)=0 \]
with coefficients $a_{ij}^k,b_{ij}$ holomorphic on $M$ and a local solution $f$ holomorphic around some $p \in M$.
For $p \in M$ the space of local solutions near $p$ is a complex vector space which we denote by $V_p^{\ast}$. 
Clearly a local solution $f \in V_p^{\ast}$ is completely determined 
by the numerical values $f(p)$ and $\partial_i f(p)$ for $i=1,\cdots,n$.
Hence the dimension of $V_p^{\ast}$ is at most equal to $(n+1)$. 

\begin{definition}
The system of $n(n+1)/2$ linear differential equations
\[ [\partial_i\partial_j+\sum_k\;a_{ij}^k(z)\partial_k+b_{ij}(z)]f(z)=0 \]
is called integrable if the dimension of $V_p$ is equal to $(n+1)$ for all $p \in M$.
\end{definition}

We assume from now on that the system is integrable.
If $\gamma:[0,1] \rightarrow M$ is a path from $p=\gamma(0)$ to $q=\gamma(1)$
then analytic continuation of local solutions along $\gamma$ gives an invertible linear map
\[ A_{\gamma}^{\ast}: V_p^{\ast} \rightarrow V_q^{\ast} \]
of local solution spaces.
Indeed the inverse of $A^{\ast}_{\gamma}$ is equal to $A^{\ast}_{{\gamma}^{-1}}$ 
with $\gamma^{-1}$ the path $\gamma$ in reversed direction.
If the path $\delta:[0,1] \rightarrow M$ has begin point $q=\delta(0)$ and end point $r=\delta(1)$
then the composition $\delta\gamma$ is defined by first traversing $\gamma$ and subsequently $\delta$.
It has begin point $p$ and end point $r$ and the product rule
\[ A^{\ast}_{\delta\gamma}=A^{\ast}_{\delta}A^{\ast}_{\gamma}: V_p^{\ast} \rightarrow V_r^{\ast} \]
is obvious. Let us write $\Pi_1(M,p)$ for the fundamental group of $M$ with base point $p$. 
The product rule on $\Pi_1(M,p)$ is derived from the above composition rule.
Hence the map $\gamma \mapsto A^{\ast}_{\gamma}$ yields a representation
\[ \rho^{\ast}: \Pi_1(M,p) \rightarrow \GL(V_p^{\ast}) \]
which is called the dual monodromy representation of the integrable system.

Let us assume that the base point $p \in M$ is fixed throughout the rest of these notes.
So we will drop the index $p$ and simply write $V^{\ast}$, $\Pi_1(M)$ and so on.
Let us denote $V=\Hom(V^{\ast},\C)$ for the vector space dual to $V^{\ast}$, and write
\[ \rho: \Pi_1(M) \rightarrow \GL(V)\;,\; \rho(\gamma)=(\rho^{\ast}(\gamma^{-1}))^{\ast}=A_{{\gamma}^{-1}} \]
for the monodromy representation.
The image $G=\rho(\Pi_1(M))$ inside $\GL(V)$ is called the monodromy group,
and we write $\Gamma=PG$ inside $\PGL(V)$ for the projective monodromy group acting on $\Pro(V)$.
Write $P:\GL(V) \twoheadrightarrow \PGL(V)$ for the natural map.
Note that $\Pi_1(M)/\Ker(P\rho) \cong \Gamma$ by the homomorphism theorem.

We define the multivalued evaluation map
\[ ev: M \dashrightarrow V\;,\;f(ev(q))=A^{\ast}_{\gamma}f(q) \]
for all $f \in V^{\ast}$. Here $\gamma: [0,1] \rightarrow$ is a path from $p$ to $q$,
and $A^{\ast}_{\gamma}f$ is the analytic continuation of $f \in V^{\ast}$ along $\gamma$ as above.
The multivaluedness of the evaluation map comes from the choice of a path $\gamma$ from $p$ to $q$.

Note that $ev(q)$ is a nonzero vector in $V$ for all $q \in M$, 
because for each $q \in M$ there exists a local solution in $V_q^{\ast}$ which does not vanish at q.
Hence the multivalued projective evaluation map 
\[ Pev: M \dashrightarrow \Pro(V)\;,\;Pev(q)=[ev(q)] \]
is well defined. 
Here we denote by $[v] \in \Pro(V)$ the line through $v \in V-\{0\}$.
Since the projective evaluation map is the main object of study 
we have chosen to denote the initial local solution space by $V^{\ast}$ 
and the dual (range of evaluation map) vector space by $V$.

In order to eliminate the multivaluedness of these maps let us denote by $\widehat{M}$
the universal $\Gamma$-covering space of $M$ related to our fixed base point $p$.
In other words, $\widehat{M}$ is equal to $\Ker(P\rho)\backslash\widetilde{M}$ 
with $\widetilde{M}$ the universal covering of $M$, and $\Gamma \cong \Pi_1(M)/\Ker(P\rho)$ the projective monodromy group. 
Then we have the commutative diagram
\[ \begin{CD}
      \widehat{M}   @>\widehat{Pev}>>       \Pro(V) \\
         @VVV                          @VVV         \\
      M             @>Pev>> \Gamma\backslash\Pro(V) \\ 
   \end{CD} \]
with the left vertical arrow an unramified $\Gamma$-covering map, 
and the group $\Gamma=\Pi_1(M)/\Ker(P\rho)$ acting as group of deck transformations on the left side of the diagram, 
and the group $\Gamma=PG$ acting as projective monodromy group on the right side of the diagram. 
The top horizontal arrow is singlevalued and removes the multivaluedness of the projective evaluation map on $M$.
But $\Gamma\backslash\Pro(V)$ is an ill defined space unless the action of $\Gamma$ 
on (the image of $\widehat{Pev}$ in) $\Pro(V)$ is properly discontinuous.
Recall that the multivalued Wronskian (with $\partial_0f_j=f_j$)
\[ W(f_0,\cdots,f_n)=\det(\partial_i f_j)_{0 \leq i,j \leq n} \]
is nowhere vanishing if and only if $f_0,f_1,\cdots,f_n$ is a basis of the local solution space $V^{\ast}$.
Note that the Wronskian is a solution of the first order system
\[ [\partial_j+\sum_{i=1}^n a_{ij}^i(z)]W(f_0,\cdots,f_n)(z)=0 \]
in local coordinates as before. 
Hence once the Wronskian is nonzero somewhere it remains nonzero everywhere.
A basis $f_0,\cdots,f_n$ of $V^{\ast}$ identifies $\Pro(V)$ with $\Pro^n(\C)$,
and identifies 
\[ Pev(q) \simeq [f_0(q):f_1(q): \cdots :f_n(q)] : \widehat{M} \rightarrow \Pro^n(\C) \]
with the projective evaluation map. This map is locally biholomorphic since (say $f_0 \neq 0$)
\[ W(f_0,\cdots,f_n)=f_0^{n+1}\det(\partial_i(f_j/f_0))_{1 \leq i,j \leq n} \]
is essentially the Jacobian of the projective evaluation map in the affine chart $\{f_0 \neq 0\}$.
The conclusion is that $M$ is locally modelled on the projective space $\Pro(V)$ of dimension $n$. 
We say that $M$ is equipped with a projective structure,
and the projective evaluation map is also called the Schwarz map of the projective structure.

\begin{definition}
A projective structure on $M$ with Schwarz map 
\[ \widehat{Pev}: \widehat{M} \rightarrow \Pro(V) \]
is called an elliptic structure if $V$ has a Hermitian form $\langle \cdot,\cdot \rangle$ of Euclidean signature
that is invariant under the monodromy group $G$.
\end{definition}

In other words $M$ becomes a K\"{a}hler manifold locally modelled on 
the projective space $\Pro(V)$ equipped with the Fubini-Study metric. 
By construction the Schwarz map $\widehat{Pev}: \widehat{M} \rightarrow \Pro(V)$ becomes a local isometry.

Suppose $V$ has a Hermitian form $\langle \cdot,\cdot \rangle$ of Lorentzian signature $(n,1)$ then
\[ \B(V) = \{[v]\in \Pro(V);\langle v,v \rangle <0 \} \]
is called the complex hyperbolic ball of dimension $n$. 
The complex ball is endowed with its natural hyperbolic metric.

\begin{definition}
A projective structure om $M$ with Schwarz map 
\[ \widehat{Pev}: \widehat{M} \rightarrow \Pro(V) \]
is called a hyperbolic structure if $V$ has a Hermitian form $\langle \cdot,\cdot \rangle$ of signature $(n,1)$
that is invariant under the monodromy group $G$. Moreover we also require that 
\[ \widehat{Pev}: \widehat{M} \rightarrow \B(V) \]
so the image of the Schwarz map should be contained inside the ball $\B(V)$.
\end{definition}

In other words $M$ becomes a K\"{a}hler manifold, which is locally modelled on the complex hyperbolic ball $\B(V)$.
By construction the Schwarz map $\widehat{Pev}: \widehat{M} \rightarrow \B(V)$ is a local isometry.
An elliptic or hyperbolic structure is called a geometric structure,
and the Schwarz map is also called the developing map of the geometric structure.
The intuitive idea of the developing map is the unrolling 
of the geometric manifold $M$ over $\Pro(V)$ or $\B(V)$ respectively.
If the vector space $V$ with the Hermitian form of Euclidean or 
Lorentzian signature is given we simply write $\Pro=\Pro(V)$ and $\B=\B(V)$.

\begin{example}
Suppose $M$ has a hyperbolic structure, for which the hyperbolic metric on $M$ is complete.
For example if $M$ is a compact space the metric is always complete. 
Never mind that this contradicts our (avoidable) assumption that $M \subset \C^n$ is Zariski open.
Then the developing map 
\[ \widehat{Pev}: \widehat{M} \rightarrow \B \]
becomes an unramified covering map, hence is a biholomorpic isomorphism, since $\B$ is simply connected.
Therefore we get a commutative diagram
\[ \begin{CD}
      \widehat{M} @>\widehat{Pev}>>       \B  \\
         @VVV                           @VVV  \\
      M           @>Pev>>  \Gamma\backslash\B \\ 
   \end{CD} \]
with $\Gamma$ a discrete torsionfree subgroup of $\Aut(\B)\cong\PU(V,\langle \cdot,\cdot \rangle)$.
Any $M=\Gamma\backslash\B$ for such $\Gamma$ arrises this way, and such $\Gamma$ exist in abundance.
\end{example}

We are interested in the general setting of a smooth projective manifold $\overline{M} = M \sqcup D$
with $D$ a divisor that is an arrangement, that is in suitable local coordinates
$D$ is a finite collection of hyperplanes. Suppose we have given on $M$ a projective structure
\[ [\partial_i\partial_j+\sum_k\;a_{ij}^k(z)\partial_k+b_{ij}(z)]f(z)=0 \]
which is regular singular along $D$. This means that local solutions have moderate growth along $D$.
Assume that the projective structure on $M$ is in fact a hyperbolic structure.
Let $G$ be the monodromy group and $\Gamma=PG$ the projective monodromy group.

In our examples it seldom happens that the hyperbolic metric on $\widehat{M}$ is complete.
However under very restrictive conditions it is possible to construct a new compactification $M \hookrightarrow M^+$
with an associated partial compactifications $\widehat{M} \hookrightarrow \widehat{M}^+$, 
and a partial compactification $\B \hookrightarrow \B^+$, both with action by $\Gamma$, 
and a continuous equivariant extension 
$\widehat{Pev}^+:\widehat{M}^+ \longrightarrow \B^+$ such that the diagram
\[ \begin{CD}
      \widehat{M}   @>>> \widehat{M}^+   @>\widehat{Pev}^+>>   \B^+    \\
         @VVV               @VVV                                @VVV   \\
      M             @>>>  M^+            @>Pev^+>>\Gamma\backslash\B^+ \\ 
   \end{CD} \]
commutes. 
The left vertical arrow is an unramified covering, while the other two vertical arrows should be ramified coverings.
The crucial property we want is that the extended developing map
\[ \widehat{Pev}^+: \widehat{M}^+ \rightarrow \B^+ \]
is a locally biholomorphic covering map over the full ball $\B$, and hence is a biholomorphic isomorphism.
The desired compactification $M^+$ of $M$ should be obtained from the divisor $D$ in $\overline{M}$ 
by a birational Cremona type transformation (blow up and blow down process).
By a careful analysis of the local exponents of our given integrable system corresponding to the hyperbolic structure
one can see which of the strata of the intersection lattice of $D$ need to be blown up and in which order.
However the blow down process can only be achieved in an indirect way from Stein factorization in the diagram
\[ \begin{CD}
      \widehat{M}   @>>> \widehat{M}^+  @>>> \widehat{M}^+_{\mathfrak{St}}  @>\widehat{Pev}^+_{\mathfrak{St}}>>   \B^+    \\
         @VVV               @VVV                                @VVV                                         @VVV         \\
      M             @>>>  M^+           @>>>  M^+_{\mathfrak{St}}           @>Pev^+_{\mathfrak{St}}>>\Gamma\backslash\B^+ \\ 
   \end{CD} \]
Again we want the Stein factorization
\[ \widehat{Pev}^+_{\mathfrak{St}}: \widehat{M}^+_{\mathfrak{St}} \longrightarrow \B^+  \]
to be a locally biholomorphic covering over the full ball $\B$, and hence a biholomorphic isomorphism.

It is absolutely necessary to understand the one dimensional case of the Euler-Gauss hypergeometric equation,
before dealing with the delicacies of the higher dimensions.
  
\section{The Euler-Gauss hypergeometric equation}

Let $M$ be equal to $\Pro-\{0,1,\infty\}$ with variable $z$, $\theta=z\partial$ and $\Pro=\C\sqcup\{\infty\}$ the Riemann sphere.
The Euler-Gauss hypergeometric equation is the second order linear differential equation
\[ [\theta(\theta+\gamma-1)-z(\theta+\alpha)(\theta+\beta)]f=0 \]
with $\alpha,\beta,\gamma$ three complex parameters.
We shall be interested in the case that $\alpha,\beta,\gamma \in \R$ and even that $\alpha,\beta,\gamma \in \Q$.
It has regular singular points at $z=0,1$ and $\infty$ with local exponents given by the Riemann scheme
\begin{center}
\begin{tabular}{|l|l|l|} \hline
$0$ & $1$ & $\infty$ \\ \hline
$0$ & $0$ & $\alpha$ \\ \hline
$1 - \gamma$ & $\gamma - (\alpha+\beta)$ & $\beta$  \\ \hline
\end{tabular}
\end{center}
The first line contains the three singular points and the next two lines give the local exponents at these points. 
So the exponent differences at $0$, $1$ and $\infty$ are 
$\kappa = 1-\gamma$, $ \lambda = \gamma-(\alpha+\beta)$ and $\mu = \beta-\alpha$ respectively.
We shall assume that $0 \leq \kappa,\lambda,\mu$ and $\kappa+\lambda,\kappa+\mu,\lambda+\mu \leq 1$, 
which can always be arranged after shifting $\alpha,\beta,\gamma$ by integers and 
performing sign changes on the differences $\kappa,\lambda,\mu$. Such a parameter set is called reduced.

Let $\Pi=\Pi_1(M,p)$ be the fundamental group of $M=\Pro-\{0,1,\infty\}$ with base point $p=1/2$
with three generators $\gamma_0,\gamma_1,\gamma_{\infty}$
and the single relation $\gamma_{\infty}\gamma_1\gamma_0=1$ as indicated in the following picture.

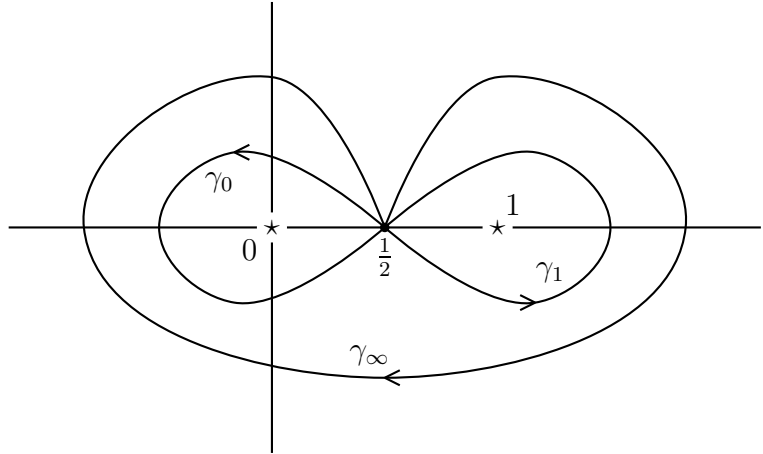
\begin{center}
\psset{unit=1mm}
\begin{pspicture}*(-50,-30)(50,30)

\psline(-50,0)(50,0)
\psline(-15,-30)(-15,30)
\pscurve(0,0)(-20,10)(-30,0)(-20,-10)(0,0)
\psline(-20,10)(-18,11)
\psline(-20,10)(-18,9)
\pscurve(0,0)(20,10)(30,0)(20,-10)(0,0)
\psline(20,-10)(18,-11)
\psline(20,-10)(18,-9)
\pscurve(0,0)(-15,20)(-40,0)(0,-20)(40,0)(15,20)(0,0)
\psline(0,-20)(2,-21)
\psline(0,-20)(2,-19)
\pscircle*[linecolor=white](-15,0){2}
\pscircle*[linecolor=white](15,0){2}

\rput(-2,-17){$\gamma_{\infty}$}
\rput(-22,6){$\gamma_0$}
\rput(22,-6){$\gamma_1$}
\rput(-18,-3){$0$}
\rput(17,3){$1$}
\rput(0,-4){$\frac{1}{2}$}
\rput(-15,0){$\star$}
\rput(15,0){$\star$}
\psdot(0,0)

\end{pspicture}
\end{center}

Now let us pick two linearly independent solutions $f_1,f_2$ on the upper half plane $\mathbb{H}$,
and consider the projective evaluation map
\[ Pev: \UH \rightarrow \Pro\;,\;Pev(z)=f_1(z)/f_2(z)  \] 
which is also called the Schwarz map.
Because of the ambiguity of the base choice $f_1,f_2$ the Schwarz map 
is only canonical up to an action of $\mbox{Aut}(\mathbb{P})$.
We claim that the the Schwarz map $Pev$ maps the upper half plane $\UH$ conformally
onto the interior of a triangle with sides circular arcs, and with angles 
$\kappa \pi$, $\lambda \pi$ and $\mu \pi$ at the vertices $Pev(0)$, $Pev(1)$ and $Pev(\infty)$ respectively.
This circular triangle is called the Schwarz triangle of the hypergeometric equation.

For example, for the boundary interval $(0,1)$ let us choose the solutions $f_1,f_2$ to be real on $(0,1)$.
This is possible since the hypergeometric equation is a real differential equation,
since the parameters $\alpha,\beta,\gamma$ were assumed to be real numbers.
In that case the image of the interval $(0,1)$ under the Schwarz map is a real interval.
For a general choice of $f_1,f_2$ the image of $(0,1)$ is the transform under an element of $\mbox{Aut}(\mathbb{P})$,
so a fractional linear transformation, of a real interval, and therefore equal to a real interval or a circular arc.

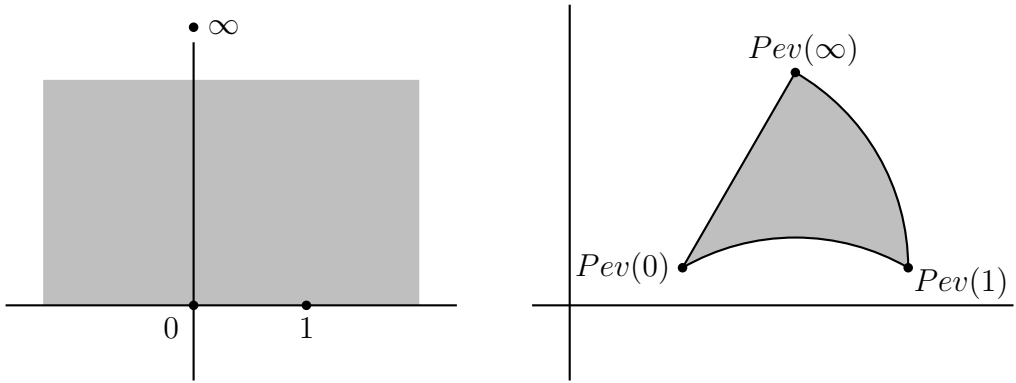
\begin{center}
\psset{unit=1mm}
\begin{pspicture}*(-65,-20)(69,30)

\pspolygon*[linecolor=lightgray](-60,-10)(-10,-10)(-10,20)(-60,20)
\pswedge*[linecolor=lightgray](25,-5){30}{0}{60}
\pswedge*[linecolor=white](40,-30.98){30}{60}{120}

\psarc(25,-5){30}{0}{60}
\psarc(40,-30.98){30}{60}{120}

\psline(25,-5)(40,20.98)
\psline(-65,-10)(-5,-10)
\psline(-40,-20)(-40,25)

\psdot(-40,-10)
\psdot(-25,-10)
\psdot(-40,27)
\psdot(25,-5)
\psdot(55,-5)
\psdot(40,20.98)

\rput(-43,-13){0}
\rput(-25,-13){1}
\rput(-36,27){$\infty$}

\psline(5,-10)(69,-10)
\psline(10,-20)(10,30)

\rput(17,-5){$Pev(0)$}
\rput(62,-7){$Pev(1)$}
\rput(41,24){$Pev(\infty)$}

\end{pspicture}
\end{center}

The angles of the Schwarz triangle at the vertices $Pev(0)$, $Pev(1)$ and $Pev(\infty)$
are equal to $\kappa \pi$, $\lambda \pi$ and $\mu \pi$ respectively.
For example, near the origin $0$ let us choose the solutions $f_1,f_2$ of the form
\[ f_1(z)=z^{1-\gamma}(1+ \cdots)\;,\; f_2(z)=(1+ \cdots) \]
which in turn implies that 
\[ f_1(z)/f_2(z)=z^{\kappa}(1+ \cdots) \]
which indeed gives an angle $\kappa \pi$ at the vertex $Pev(0)$ of the Schwarz triangle.
For a general choice of $f_1,f_2$ this angle $\kappa \pi$ is conserved by some fractional linear transformation.

By continuity we can extend the Schwarz map
\[ Pev: \mathbb{H} \sqcup (-\infty,0) \sqcup (0,1) \sqcup (1,\infty) \rightarrow \mathbb{P} \]
with image the Schwarz triangle minus its vertices.
The key step in the argument of Schwarz is the beautiful insight 
that the analytic continuation of $Pev$ is given by the reflection principle.
Indeed, there are three possibilities for analytic continuation 
from the upper half plane $\mathbb{H}$ to the lower half plane $-\mathbb{H}$,
namely through the intervals $(-\infty,0)$, $(0,1)$ and $(1,\infty)$.
The analytic continuation of the Schwarz map is obtained 
by reflecting the Schwarz triangle in the corresponding sides.
Now we can iterate the above construction with the new triangle,
which allows one to understand the full analytic continuation of the Schwarz map,
step by step reflecting in sides of circular triangles.
The domain of this full analytic continuation is the universal covering space 
$\widetilde{M}$ of $M=\Pro-\{0,1,\infty\}$, and we write
\[ \widetilde{Pev}: \widetilde{M} \rightarrow \mathbb{P} \]
for the analytic continuation of the Schwarz map. 

The range of this map can get messy, as the triangles start overlapping.
However in case the Schwarz triangle is dihedral, which means that 
\[ \kappa=1/k \;,\; \lambda=1/l \;,\; \mu=1/m \] 
for some integers $k,l,m \geq 1$,
we do get a regular tesselation by mirror images of the Schwarz triangle.
These requirements are called the Schwarz integrality conditions.
The range $\mathbb{G}$ of this tesselation is equal to 
\[ \mathbb{G}=\mathbb{P},\; \mathbb{C},\; \mathbb{D} \]
upto an action of $\mbox{Aut}(\mathbb{P})$, depending on whether the angle sum of the Schwarz triangle
\[ (\kappa + \lambda + \mu)\pi=(1/k+1/l+1/m)\pi \] 
is greater than $\pi$, equal to $\pi$, or smaller than $\pi$ respectively.
Here 
\[ \mathbb{D}=\{w \in \mathbb{C};\abs{w}<1\} \] 
denotes the unit disc.
In this last case the disc $\mathbb{D}$ is bounded by a circle (Klein's Orthogonalkreis)
which is orthogonal to the three circles bounding the Schwarz triangle. 
This leads one to spherical, Euclidean and hyperbolic planar geometry respectively.
Note that in all three cases $\mathbb{G}$ is simply connected.
These results go back to Riemann \cite{Riemann}, Schwarz \cite{Schwarz} and Klein \cite{Klein} 
in the golden $19^{th}$ century of German mathematics. 

Under the Schwarz conditions we shall now rephraze the above results in the language used in the first section.
The fundamental group $\Pi$ of $M=\mathbb{P}-\{0,1,\infty\}$ with base point $z_0=1/2$
has three generators $\gamma_0,\gamma_1,\gamma_{\infty}$ with the single relation $\gamma_{\infty}\gamma_1\gamma_0=1$.
The projective monodromy group $\Gamma=\Gamma(k,l,m)$ is isomorphic to a quotient of $\Pi/\Pi(k,l,m)$ with
$\Pi(k,l,m)$ be the normal subgroup of $\Pi$ generated by $\gamma_0^k$, $\gamma_1^l$ and $\gamma_{\infty}^m$.

It is clear that the Schwarz map factors through the intermediate covering 
$\widehat{M}=\Pi(k,l,m) \backslash \widetilde{M}$, and so the Schwarz map
\[ \widehat{Pev}: \widehat{M} \rightarrow \mathbb{G} \]
is locally biholomorphic with image the complement in $\mathbb{G}$ of all vertices of the tesselating triangles.
At this point we can fill in the vertices and the points above the vertices of the tesselating triangles,
resulting in a commutative diagram 
\[ \begin{CD}
      \widehat{M}=\Pi(k,l,m) \backslash \widetilde{M} @>>> \widehat{M}^+      @>\widehat{Pev}^+>> \mathbb{G}      \\
         @VVV               @VVV                     @VVV                     \\
      M=\Pro-\{0,1,\infty\}     @>>>  M^+=\Pro     @>Pev^+>>\Gamma\backslash\mathbb{G} \\ 
   \end{CD} \]
The left vertical arrow is an unramified $\Gamma$-covering, 
while the middle vertical arrow is a ramified covering, 
branched of orders $k$, $l$ and $m$ above the points $0$, $1$ and $\infty$ respectively.
The extended Schwarz map 
\[ \widehat{Pev}^+: \widehat{M}^+ \rightarrow \mathbb{G} \]
is a locally biholomorphic covering map. 
Since $\widehat{M}^+$ is connected and $\mathbb{G}$ simply connected such a covering is a bijection, 
and we conclude that the extended Schwarz map is biholomorphic from $\widehat{M}^+$ onto $\mathbb{G}$.

In other words the projective monodromy group $\Gamma(k,l,m) \cong \Pi/\Pi(k,l,m)$ acts properly discontinuously
on $\mathbb{G}$ with quotient $\Gamma(k,l,m) \backslash \mathbb{G} \cong \Pro$. 
This quotient map is ramified above the images under $Pev^+$ of $0,1,\infty$ of orders $k,l,m$ respectively.
The projective monodromy group $\Gamma(k,l,m)$ is called the Schwarz triangle group.

The group $W(k,l,m)$ generated by the reflections in the sides of the Schwartz triangle is called the Coxeter triangle group.
The Schwartz triangle group is the index two subgroup of the Coxeter triangle group,
consisting of even products of reflections in the sides of the Schwartz triangle.
The Coxeter triangle group $W(k,l,m)$ acts on $\mathbb{G}$ by holomorphic
and antiholomorphic transformations. 

\section{Binary lift of the hypergeometric equation}

Consider the map $w = (\frac12-\frac14(z+1/z)) = -(z^{\frac12}-z^{-\frac12})^2/4$ from $\Pro$ to $\Pro$.
It is a degree $2$ ramified covering as indicated in the following picture
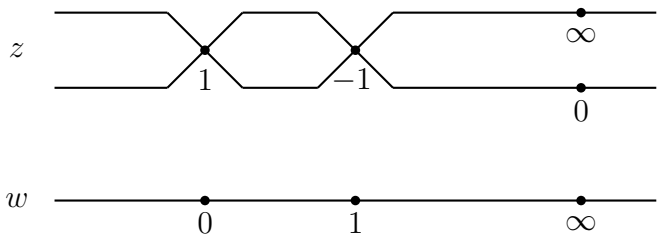
\begin{center}
\psset{unit=1mm}
\begin{pspicture}*(-40,-6)(50,30)

\psline(-30,0)(50,0)

\psline(-30,15)(-15,15)
\psline(-15,15)(-5,25)
\psline(-5,25)(5,25)
\psline(5,25)(15,15)
\psline(15,15)(50,15)
\psline(-30,25)(-15,25)
\psline(-15,25)(-5,15)
\psline(-5,15)(5,15)
\psline(5,15)(15,25)
\psline(15,25)(50,25)

\psdot(-10,0)
\psdot(10,0)
\psdot(40,0)
\psdot(-10,20)
\psdot(10,20)
\psdot(40,25)
\psdot(40,15)

\rput(-10,-3){$0$}
\rput(10,-3){$1$}
\rput(40,-3){$\infty$}
\rput(40,12){$0$}
\rput(9.5,16){$-1$}
\rput(-10,16){$1$}
\rput(40,22){$\infty$}
\rput(-35,20){$z$}
\rput(-35,0){$w$}

\end{pspicture}
\end{center}
with ramification points $z=1$ above $w=0$ and $z=-1$ above $w=1$.
Consider the Euler-Gauss hypergeometric equation on the projective line $\Pro$ with coordinate $w$
and with Riemann scheme
\begin{center}
\begin{tabular}{|l|l|l|} \hline
$w=0$ & $w=1$ & $w=\infty$ \\ \hline
$0$ & $0$ & $\alpha$ \\ \hline
$1 - \gamma$ & $\gamma - (\alpha+\beta)$ & $\beta$  \\ \hline
\end{tabular}
\end{center}
of regular singular points and local exponents.
The pull back under the map $w = \frac12-\frac14(z+z^{-1})$ of the hypergeometric equation 
lives on the projective line $\Pro$ with coordinate $z$ and takes the form
\[ [ \theta^2 + k_1\frac{1+z^{-1}}{1-z^{-1}}\theta + 
   2k_2\frac{1+z^{-2}}{1-z^{-2}}\theta + (\tfrac12k_1+k_2)^2-\lambda^2 ] f=0 \]
with $\theta=z\partial$ and $\partial=d/dz$ as before. It is easy to check the linear relations
\[ \alpha=\lambda+\tfrac12k_1+k_2\;,\; \beta=-\lambda+\tfrac12k_1+k_2\;,\;\gamma=\tfrac12+k_1+k_2 \]
between the two parameter sets. This pull back has four regular singular points 
$z=1,-1,0,\infty$ with Riemann scheme
\begin{center}
\begin{tabular}{|l|l|l|l|} \hline
$z=1$ & $z=-1$ & $z=0$ & $z=\infty$ \\ \hline
$0$ & $0$ & $\alpha$ & $\alpha$ \\ \hline
$2-2\gamma$ & $2\gamma-2(\alpha+\beta)$ & $\beta$ & $\beta$ \\ \hline
\end{tabular}
\end{center}
as is clear from the above ramification picture and 
the Riemann scheme of the Euler-Gauss hypergeometric equation.

The multiplicative group $\C^{\times}$ with the action of the group $S_2=\{\pm1\}$ by $z\mapsto z^{\pm1}$
together with the pull back of the hypergeometric equation has a natural generalization. 
Let $G$ be a simple complex Lie group with maximal torus $H$ and Weyl group $W$.
Instead of $\C^{\times}$ with the action of $S_2$ we consider 
the complex torus $H\cong(\C^{\times})^n$ with the action of $W$.
For simplicity we restrict ourselves to the simply laced case of type $ADE$.

\section{The root system hypergeometric equation}

Let $H\cong(\C^{\times})^n$ be a complex torus of dimension $n$ with Lie algebra $\mathfrak{h}\cong\C^n$.
For $\lambda$ in the character lattice $Q=\Hom(H,\C^{\times})$ we denote the associated
character by $H\ni h\mapsto h^{\lambda}\in\C^{\times}$. 
For $\xi$ in the Lie algebra $\mathfrak{h}$ we denote the associated translation
invariant vector field by $\theta_{\xi}$.  
The relation
\[ \theta_{\xi}(h^{\lambda})=\lambda(\xi)h^{\lambda}\;\;\forall h\in H \]
identifies $Q$ with a sublattice of $\mathfrak{h}^{\ast}=\Hom(\mathfrak{h},\C)$.

Assume that $Q$ carries an even integral scalar product $(\cdot,\cdot)$ such that 
\[ R=\{\alpha\in Q;(\alpha,\alpha)=2\} \]
is an irreducible root system of type $ADE$ of rank $n$.
Extend the scalar product from $Q$ to $\mathfrak{h}^{\ast}$ and transfer it to $\mathfrak{h}\cong\mathfrak{h}^{\ast}$.
Let $W<\Aut(\mathfrak{h}^{\ast})$ be the Weyl group generated by the
orthogonal reflections $s_{\alpha}(\lambda)=\lambda-(\lambda,\alpha)\alpha$ on $\mathfrak{h}^{\ast}$,
and transfer the action of $W$ to $\mathfrak{h}$ and $H$ as well.

Let us fix a decomposition $R=R_+\sqcup R_-$ into positive and negative roots.
The full automorphism group $W'$ of $Q$ is equal to $W'=W \rtimes D$ 
with $D$ the subgroup of $W'$ fixing the associated positive chamber.

The codimension one subtorus $H_{\alpha}=\{h\in H;h^{\alpha}=1\}$
is called the mirror associated with the root $\alpha$ in $R$.
We denote $H=H^{\circ}\sqcup(\cup H_{\alpha})$ and call $H^{\circ}$ 
the complement of the toric mirror arrangement.

\begin{theorem}
The system of n(n+1)/2 linearly independent differential equations
\[ [\theta_{\xi}\theta_{\eta}+\tfrac{1}{2}k\sum_{\alpha>0}\alpha(\xi)\alpha(\eta)\frac{1+h^{-\alpha}}{1-h^{-\alpha}}
   \theta_{\alpha^{\lor}}+ak^2(\xi,\eta)]f(h)=0 \;\;\forall \xi,\eta\in \mathfrak{h}\]
is an integrable system on the toric mirror arrangement complement $H^{\circ}$, and is invariant under $W'$.
Here $\mathfrak{h}^{\ast}\cong\mathfrak{h}$ via $(\cdot,\cdot)$ and 
$\mathfrak{h}^{\ast}\ni\alpha\cong\alpha^{\lor}\in\mathfrak{h}$.
Moreover $a=(n+1)/4,n-2,6,12,30$ for type $A_n,D_n,E_6,E_7,E_8$ respectively.
\end{theorem}

As discussed in the first section the above special hypergeometric system associated with the root system $R$
gives a projective structure on $H^{\circ}$, which is invariant under $W'$.

In the example of the root system of type $A_1$ this equation boils down to 
(take $\xi=\eta=\alpha^{\lor}/2$ with variable $z=h^{\alpha}$ and derivative $\theta=z\partial$)
\[ [\theta^2+k\frac{1+z^{-1}}{1-z^{-1}}\theta+\tfrac{1}{4}k^2]f(z)=0 \]
which is the pull back of the Euler-Gauss hypergeometric equation with the parameters $k_1=k$, $k_2=0$ and $\lambda=0$
as described in the previous section.

\begin{lemma}
The orbifold fundamental group $\Pi_1(W'\backslash H^{\circ})$ is the extended affine Artin group
\[ \widetilde{\Art}'=\widetilde{\Art}\rtimes \widetilde{D} \] 
with $\widetilde{\Art}$ the affine Artin group with the usual generators $T_0,T_1,\cdots,T_n$
and braid relations, and $\widetilde{D}$ the group of automorphisms of the fundamental alcove.
\end{lemma}

\begin{lemma}
In the monodromy representation $\rho:\Pi_1(W'\backslash H^{\circ})\rightarrow\GL_{n+1}(\C)$
the Artin generators $T_i$ map to complex reflections $t_i=\rho(T_i)$ satisfying
the quadratic Hecke relation $(t_i-1)(t_i+q)=0$ with $q=\exp(-2\pi ik)$.
\end{lemma}

\begin{lemma}
For $0<k<m$ with $m$ the hyperbolic exponent given by
\begin{gather*}
   m = \begin{cases}
         2/(n+1) & \mbox{for type}\; A_n \\
         1/(n-2) & \mbox{for type}\; D_n \\
         1/(n-3) & \mbox{for type}\; E_n
         \end{cases}
\end{gather*}
the monodromy representation is irreducible, and has an
invariant Hermitian form of Lorentz signature.
\end{lemma}

\begin{lemma}
For $0<k<m$ the image of the Schwarz map
\[ \widehat{Pev}: \widehat{W'\backslash H^{\circ}} \rightarrow \Pro^n=\Pro^n(\C) \]
is contained in the ball $\B^n=\B^n(\C)$.
\end{lemma}

The conclusion is that for each $0<k<m$ the toric mirror arrangement complement $H^{\circ}$ 
is endowed with a hyperbolic structure, which is invariant under $W'$.
Equivalently $W'\backslash H^{\circ}$ has a hyperbolic structure as an orbifold.

\section{The Schwarz conditions}

Let us assume that the connected smooth complex manifold $M$ has a smooth projective
compactification $\overline{M}=M\sqcup D$ with $D$ a divisor with normal crossings.
In local coordinates $z=(z_1,\cdots,z_n)$ the divisor $D$ is given by $z_1\cdots z_d=0$ for some $d=1,\cdots,n$.
Suppose that we have given on $M$ a hyperbolic structure with regular singularities along $D$, which takes the form
\[ [\partial_i\partial_j+\sum_k\;a_{ij}^k(z)\partial_k+b_{ij}(z)]f(z)=0 \]
in these local coordinates. 
Suppose this system has local exponents (an unordered set with with possible repetitions) $\mu_1,\cdots,\mu_{n+1}$ in $\Q^d$, 
where we write $z^{\mu}=z_1^{m_1}\cdots z_d^{m_d}$ if $\mu=\sum m_j\epsilon_j$
with $\{\epsilon_1,\cdots,\epsilon_d\}$ the standard basis of $\Z^d$.
Clearing out common denominators in a minimal way by going to a finite covering we can assume that
the local exponents $\{\mu_1,\cdots,\mu_{n+1}\}$ are in fact contained in $\Z^n$.
In turn, this approximates locally on $M$ the projective evaluation map
\[ z=(z_1,\cdots,z_n) \mapsto (z^{\mu_1}:z^{\mu_2}:\cdots:z^{\mu_{n+1}}) \]
as a rational map with possibly poles along the divisor $D$.

\begin{definition}
We say that the compactification $\overline{M}=M\sqcup D$ is well adapted to the given projective structure on $M$ 
if (on a finite covering as above) the projective evaluation map extends locally over $D$ as a rational map. 
This means that the closure and inclusion relations on the intersection lattice of $D$ 
are preserved under the projective evaluation map.
\end{definition}

In case the compactification $\overline{M}=M\sqcup D$ is well adapted to the given projective structure
the projective evaluation map might contract some of the strata, but it has a rational extension over $D$
without having to perform any further blow ups.

\begin{definition}
Assume that the compactification $\overline{M}=M\sqcup D$ is well adapted to the given projective structure on $M$.
We say that the Schwarz conditions hold, if the projective evaluation map has locally a biholomorphic extension 
over those codimension one strata of $D$, which are not contracted under the projective evaluation map. 
\end{definition} 

The above definitions should be formulated in a slightly more general way, so as to include possibly logarithmic terms as well.
Despite the fact that in most examples logarithmic terms do occur along some strata, we leave this aside for the moment.

\begin{question}\label{Schwarz conditions}
If the hyperbolic structure om $M$ has a well adapted compactification $\overline{M}=M\sqcup D$ 
for which the Schwarz conditions do hold, then the projective monodromy group $\Gamma$ is
a discrete subgroup of $\Aut(\B)$ of cofinite volume, and 
the projective evaluation map
\[ Pev: M \hookrightarrow \Gamma\backslash\B \]
is biholomorphic onto a Heegner divisor complement.
\end{question}

\section{The Schwarz conditions for root systems}

Let $\{\alpha_1,\cdots,\alpha_n\}$ be the set of simple roots in the fixed set of positive roots $R_+$.
The biholomorpic map 
\[ H \rightarrow (\C^{\times})^n\;,\;h \mapsto z=(z_1,\cdots,z_n)=(h^{-\alpha_1},\cdots,h^{-\alpha_n}) \]
is an isomorphism of Abelian groups, 
and the embedding $(\C^{\times})^n \hookrightarrow \C^n$ defines a partial compactification of $H$.
Here is a picture for the root system of type $A_2$.
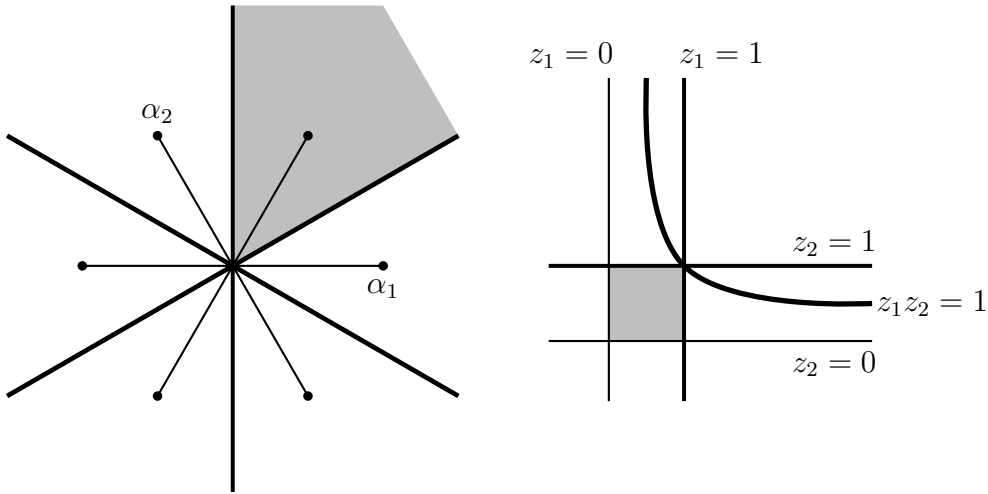
\begin{center}
\psset{unit=1mm}
\begin{pspicture}*(-60,-30)(75,40)

\psline*[linecolor=lightgray](20,-10)(30,-10)(30,0)(20,0)(20,-10)
\psline*[linecolor=lightgray](-30,0)(-30,34.64)(-10,34.64)(0,17.32)

\psline(-50,0)(-10,0)
\psline(-40,17.32)(-20,-17.32)
\psline(-40,-17.32)(-20,17.32)
\psline[linewidth=1.5pt](-30,-34.64)(-30,34.64)
\psline[linewidth=1.8pt](-60,-17.32)(0,17.32)
\psline[linewidth=1.8pt](-60,17.32)(0,-17.32)

\psline(12,-10)(55,-10)
\psline[linewidth=1.5pt](12,0)(55,0)
\psline(20,-18)(20,25)
\psline[linewidth=1.5pt](30,-18)(30,25)

\psdot(-10,0)
\psdot(-50,0)
\psdot(-40,17.32)
\psdot(-20,17.32)
\psdot(-40,-17.32)
\psdot(-20,-17.32)

\rput(-10,-3){$\alpha_1$}
\rput(-40,20.3){$\alpha_2$}
\rput(50,-13){$z_2=0$}
\rput(50,3){$z_2=1$}
\rput(63,-5){$z_1z_2=1$}
\rput(15,28){$z_1=0$}
\rput(35,28){$z_1=1$}

\pscurve[linewidth=2pt](25,25)(30,0)(55,-5)

\end{pspicture}
\end{center}
Using the action of the extended Weyl group $W'$ on the torus $H$ this partial compactification extends to 
a smooth full compactification $H \hookrightarrow \overline{H}$ by a boundary divisor with normal crossings.
These boundary divisors are called the toric strata.

The intersection of the mirror arrangement $\cup H_{\alpha}$ with this toric boundary divisor is transversal.
There is a natural way to resolve the mirror arrangement to a divisor with normal crossings.
If we are given a connected component of an intersection of mirrors, then the set of those roots, 
whose mirror contains the component, forms a root subsystem of $R$. 
We say that the connected component is irreducible if the associated root subsystem is irreducible.
We start by blowing up the irreducible connected components of dimension zero,
and consider the strict transform of the mirror arrangement on this blow up.
Subsequently we blow up the irreducible connected components of dimension one,
and consider the strict transform of the mirror arrangement on this blow up.
Repeating this construction we arrive at a smooth compactification
\[ H^{\circ} \hookrightarrow (H^{\circ})^+ \]
with a boundary divisor with normal crossings. Remark that the construction
is equivariant for the action of the extended Weyl group $W'$.

The special hypergeometric system described in the previous section has regular singular points at infinity,
both for the toric strata and the blown up mirror arrangement. 
The following results are obtained by explicit calculations.

\begin{lemma}
The Schwarz conditions for the toric strata are given for type $A_n$ by
\[ (n-1)k/2 \in 1/\N \;, \]
and for type $D_n$ or $E_n$ by
\[ dk \in 1/\N \;. \]
Here $d$ is the length in the Coxeter diagram from an extremal node to the triple node,
so $d=1,n-3$ for type $D_n$ and $d=1,2,n-4$ for type $E_n$.
\end{lemma}

\begin{lemma}
The Schwarz conditions for the mirror strata and near the identity element are given by (with $h$ the Coxeter number)
\[ (1-2k)/2 \in 1/\N \; \mbox{if} >0 \;,\; (hk-1)/2 \in 1/\N \; \mbox{if} >0 \;, \]
and near the other special point of type $A_7$ in $E_7$ by
\[ (1-2k)/2 \in 1/\N \; \mbox{if} >0 \;,\; (8k-1)/2 \in 1/\N \; \mbox{if} >0 \;, \]
and near the other special points of type $A_8$ and $D_8$ in $E_8$ by
\begin{gather*}
   (1-2k)/2 \in 1/\N \; \mbox{if} >0 \;,\; (9k-1) \in 1/\N \; \mbox{if} >0 \\
   (1-2k)/2 \in 1/\N \; \mbox{if} >0 \;,\; (14k-1)/2 \in 1/\N \; \mbox{if} >0  
\end{gather*}
respectively. Here the phrase $x\in 1/\N$ if $>0$ means if $x>0$ then $x\in 1/\N$.
\end{lemma}

\begin{corollary}
Say $(1-2k)/2=1/p$ with $p\in\N$, $p\geq3$ so the monodromy around a mirror is an order $p$ reflection.
For the root system $R$ of rank at least $2$ we find the solutions
\begin{equation*}
\begin{split}
   p=\;\;3 \;:\;&  A_2,A_3,A_4,A_7,\;  D_4,D_5,D_6,\; E_6,E_7  \\
   p=\;\;4 \;:\;&  A_2,A_3,A_5,\;     D_4,D_5,\;    E_6        \\
   p=\;\;6 \;:\;&  A_2,A_3,A_4,A_5,\; D_4                      \\
   p= 10   \;:\;& A_2     
\end{split}                    
\end{equation*}
to the above Schwarz conditions.
\end{corollary}

For all cases in the above table Question{\;\ref{Schwarz conditions}} 
says that the hyperbolic structure on $W'\backslash H^{\circ}$ is in fact 
obtained as the complement of a Heegner divisor in a ball quotient $\Gamma \backslash \B$. 

\section{A geometric interpretation}

We look for a commutative diagram
\[ \begin{CD}
   W'\backslash H^{\circ}  @>Pev>>  \Gamma'\backslash\B \\ 
         @VVV                           @VVV  \\
      \mathcal{M}          @>Per>>  \Gamma\backslash\B \\
   \end{CD} \]
with $\mathcal{M}$ a suitable moduli space and $Per$ a suitable period map.
The left vertical arrow should be a suitable finite covering map
corresponding to the finite index subgroup $\Gamma'$ of $\Gamma$.
The following theorem is due to Deligne and Mostow \cite{DM},\cite{Mostow},\cite{Thurston}, 
based on earlier work of Picard and Terada.

\begin{theorem}
Given $n+3$ distinct points $z_0,\cdots,z_{n+2}$ on the projective line $\Pro$
and	rational numbers $0<\mu_0,\cdots,\mu_{n+2}<1$ with $\sum \mu_i=2$.
Write $\mu_i=m_i/m$ as quotient of natural numbers with $\gcd(m,m_0,\cdots,m_{n+2})=1$. Consider the curve
\[ C(z):\;\; y^m=\prod (x-z_i)^{m_i} \]
with holomorphic differential $\omega=dx/y$. The periods 
\[ \int_{z_i}^{z_{i+1}} \omega \]
are solutions of the Lauricella $F_D$ hypergeometric equation.
Let $S_{\mu}$ be the subgroup of the symmetric group $S_{n+3}$ fixing $\vect{\mu}=(\mu_0,\cdots,\mu_{n+2})$.

If for all $i \neq j$ the Schwarz conditions
\[ \mu_i+\mu_j<1 \Longrightarrow (1-\mu_i-\mu_j) \in 
   \left\{ \begin{array}{ll}
           1/\N & \mbox{if $\mu_i\neq\mu_j$} \\
           2/\N & \mbox{if $\mu_i=\mu_j$}
           \end{array} 
   \right. \]
do hold then the period map 
\[ \mathcal{M}=S_{\mu}\backslash\mathcal{M}_{0,n+3} \hookrightarrow \Gamma\backslash\B \]
is an embedding as a Heegner divisor complement. 
\end{theorem}

This theorem gives a geometric interpretation for the hypergeometric system associated with the root system of type $A_n$.
Let $z_0=0,z_{n+2}=\infty$ and $z_1,\cdots,z_{n+1}\in \C^{\times}$ with $z_1 \cdots z_{n+1}=1$
and take these numbers modulo the action of the cyclic group $C_{n+1}$ of order $n+1$.
If we take $\mu_1=\cdots=\mu_{n+1}=k$ and $\mu_0=\mu_{n+2}=(1-(n+1)k/2)$ then
$W'=S_{n+1} \times S_2<S_{\mu}$ is the extended Weyl group 
and $W'\backslash H^{\circ} \rightarrow \mathcal{M}=S_{\mu}\backslash \mathcal{M}_{0,n+3}$.
The Schwarz conditions in the theorem of Deligne and Mostow (for $W'$ rather than $S_{\mu}$)
\begin{gather*}
   (1-\mu_0-\mu_1)=(n-1)k/2 \in 1/\N \\
   (1-\mu_1-\mu_{n+1})/2=(1-2k)/2 \in 1/\N \\
   (1-\mu_0-\mu_{n+2})/2=((n+1)k-1)/2 \in 1/\N
\end{gather*}
coincide with the Schwarz conditions for our special hypergeometric system
associated with the root system $A_n$ along the toric strata, along the mirrors, 
and near the identity element respectively.

Remark that $W'=S_{\mu}$ unless $k=(1-(n+1)k/2) \Leftrightarrow k=2/(n+3)$, which happens for $(p,n)=(4,5),(6,3),(10,2)$.
From the root system perspective there is additional hidden symmetry for these cases.

Let us now discuss the situation for the root system $R$ of type $E_6$ with $p=4$, 
and its relation to moduli of quartic curves.
The first step is the following result of Looijenga \cite{Looijenga}

\begin{theorem}
There is a biholomorphic isomorphism from the moduli space $\mathcal{M}'$ 
of smooth quartic curves marked with a bitangent onto
the quotient space $W'\backslash H^{\circ}$ of type $E_6$.
\end{theorem}

The second step is a result of Kondo \cite{Kondo}, 
which in turn is based on the Torelli theorem for quartic surfaces of Shafarevic and Piatetskii-Shapiro.

\begin{theorem}
There is a period isomorphism from the moduli space $\mathcal{M}$ of smooth quartic curves 
onto a Heegner divisor complement in a ball quotient $\Gamma\backslash\B$ given as follows. 
To a smooth quartic curve one associates the smooth quartic surface ramified of order $4$ over $\Pro^2$
with ramification locus the given smooth quartic curve, 
and subsequently applies the Torelli theorem for quartic surfaces.
\end{theorem}

The special hypergeometric system associated with the root system $R$ 
of type $E_6$ with $k=1/4$ is part of the following commutative diagram
\[ \begin{CD}
    \mathcal{M}' \cong W'\backslash H^{\circ}  @>Pev>>  \Gamma'\backslash\B \\ 
         @VVV                                  @VVV                         \\
      \mathcal{M}                              @>Per>>  \Gamma\backslash\B  \\
   \end{CD} \]
The left vertical arrow is the covering map of degree $28$
from the moduli space $\mathcal{M}'$ of smooth quartic curves marked with a bitangent 
onto the moduli space $\mathcal{M}$ of smooth quartic curves, by forgetting the bitangent.
The bottom horizontal arrow is the period map as obtained by Kondo. 
Again we have the phenomenon of hidden symmetry, which in this case can
be explained by looking at the moduli space $\mathcal{M}_{3,1}$ in connection
with the space $W\backslash H^{\circ}$ for the root system $R$ of type $E_7$, 
see the paper of Looijenga \cite{Looijenga}.

\noindent
Gert Heckman, Radboud University Nijmegen, P.O. Box 9010,
\newline 
6500 GL Nijmegen, The Netherlands (E-mail: g.heckman@math.ru.nl)
\newline
Eduard Looijenga, University of Utrecht, P.O. Box 80.010
\newline 
3508 TA Utrecht, The Netherlands (E-mail: looijeng@math.uu.nl)

\end{document}